\documentclass[journal]{IEEEtran}

\ifCLASSINFOpdf
  \usepackage[pdftex]{graphicx}
\else
\fi

\usepackage[justification==center]{caption}

\usepackage{mathtools}
\usepackage{amsmath}

\ifCLASSOPTIONcompsoc
 \usepackage[caption=false,font=normalsize,labelfont=sf,textfont=sf]{subfig}
\else
 \usepackage[caption=false,font=footnotesize]{subfig}
\fi

\usepackage{xspace,amssymb,epsfig,syntonly}
\usepackage{epsfig,amsmath,color}
\usepackage{xspace,syntonly,empheq}
\usepackage{url}
\usepackage{bbm}
\usepackage{mathtools}
\usepackage{cite}
\usepackage{graphicx}
\usepackage{caption}
\usepackage{algorithm,algorithmic}
\usepackage{verbatim}
\usepackage{amssymb}
\usepackage{amsthm}
\usepackage{circuitikz}
\usepackage{enumerate}

\usepackage{multirow}

\newcommand{\todo}[1]{{\color{red} TODO: #1}}

\iftrue
  \newcommand{\excludeforsubmit}[1]{}
\else
  \newcommand{\excludeforsubmit}[1]{#1}
\fi

\usepackage{amsfonts}

\usepackage{color}

\newcommand{\rui}[1]{{\color{magenta} (Rui:  #1)}}

\newcommand{\yc}[1]{{\color{blue} (YC says:  #1)}}
\newcommand{\hl}[1]{{\color{green} (Highlighted:  #1)}}

\newcommand{\slack}{1}
\newcommand{\Y}{\bY}
\newcommand{\vv}{\bv}
\newcommand{\vvm}{\bv_{-\slack}}
\newcommand{\sinj}{\bs}
\newcommand{\sinjp}[1]{\bs_{#1}}
\newcommand{\sinjm}{\sinjp{-\slack}}
\newcommand{\curr}{\bi}
\newcommand{\buses}{\cB}
\newcommand{\nbus}{|\cB|}
\newcommand{\lines}{\cL}
\newcommand{\nlines}{|\cL|}
\newcommand{\zloadv}{{\bf w}}

\newcommand{\M}{\boldsymbol{M}}
\newcommand{\X}{\boldsymbol{X}}
\newcommand{\known}{\Psi}
\newcommand{\mcrows}{n_1}
\newcommand{\mccols}{n_2}
\newcommand{\ohmsparam}{\boldsymbol{\xi}}
\newcommand{\vparam}{\boldsymbol{\tau}}
\newcommand{\vmagpar}{\boldsymbol{\gamma}}

\newcommand{\utwi}[1]{\mbox{\boldmath $#1$}}

\newcommand{\diag}{{\textrm{diag}}}

\newcommand{\reals}{\mathbb{R}}
\newcommand{\comps}{\mathbb{C}}

\newcommand{\cL}{{\cal{L}}}

\newcommand{\cB}{{\cal B}}

\newcommand{\bg}{{\bf g}}
\newcommand{\bh}{{\bf h}}

\newcommand{\bs}{{\bf s}}
\newcommand{\bt}{{\boldsymbol{t}}}
\newcommand{\bx}{{\bf x}}

\newcommand{\bv}{{\bf v}}
\newcommand{\bi}{{\bf i}}
\newcommand{\bz}{{\bf z}}
\newcommand{\by}{{\bf y}}
\newcommand{\bA}{{\bf A}}

\newcommand{\bC}{{\bf C}}

\newcommand{\bJ}{{\bf J}}

\newcommand{\bX}{{\bf X}}

\newcommand{\bY}{{\bf Y}}

\newcommand{\bepsilon}{{\utwi{\epsilon}}}

\newcommand{\sfT}{\textsf{T}}

\usepackage[usestackEOL]{stackengine}

\DeclareMathOperator*{\minimize}{minimize}
\DeclareMathOperator*{\subjectto}{subject\;to}

\usepackage{booktabs}

\newcommand\scalemath[2]{\scalebox{#1}{\mbox{\ensuremath{\displaystyle #2}}}}
\usepackage{moresize}

\begin{document}
%
\title{Matrix Completion for Low-Observability\\ Voltage Estimation}
%
%
%

\author{Priya~L.~Donti,~\IEEEmembership{Student Member,~IEEE,}
        Yajing~Liu,~\IEEEmembership{Member,~IEEE,}
        Andreas~J.~Schmitt,~\IEEEmembership{Student Member,~IEEE,\\}
        Andrey~Bernstein,~\IEEEmembership{Member,~IEEE,}
        Rui~Yang,~\IEEEmembership{Member,~IEEE,}
        and Yingchen~Zhang,~\IEEEmembership{Senior Member,~IEEE}

\thanks{Financial support from the Department of Energy Computational Science Graduate Fellowship under Grant No. DE-FG02-97ER25308 is gratefully acknowledged. This work was authored in part by the National Renewable Energy Laboratory, managed and operated by Alliance for Sustainable Energy, LLC, for the U.S. Department of Energy (DOE) under Contract No. DE-AC36-08GO28308. Funding provided by the U.S. Department of Energy Office of Energy Efficiency and Renewable Energy Solar Energy Technologies Office. The views expressed in the article do not necessarily represent the views of the DOE or the U.S. Government. The U.S. Government retains and the publisher, by accepting the article for publication, acknowledges that the U.S. Government retains a nonexclusive, paid-up, irrevocable, worldwide license to publish or reproduce the published form of this work, or allow others to do so, for U.S. Government purposes. 
}

}

\maketitle

\begin{abstract}
With the rising penetration of distributed energy resources, distribution system control and enabling techniques such as state estimation have become essential to distribution system operation. 
However, traditional state estimation techniques have difficulty coping with the \emph{low-observability} conditions often present on the distribution system due to the paucity of sensors and heterogeneity of measurements.
To address these limitations, we propose a distribution system state estimation algorithm that employs matrix completion (a tool for estimating missing values in low-rank matrices) augmented with noise-resilient power flow constraints.
This method operates under low-observability conditions where standard least-squares-based methods cannot operate, and flexibly incorporates any network quantities measured in the field.
We empirically evaluate our method on the IEEE 33- and 123-bus test systems, and find that it provides near-perfect state estimation performance (within 1\% mean absolute percent error) across many low-observability data availability regimes.
\end{abstract}


%
\IEEEpeerreviewmaketitle

\section{Introduction}
\label{sec:intro}

\IEEEPARstart{S}{tate} estimation is one of the most critical inference tasks in power systems. 
Classically, it entails estimating voltage phasors at all buses in a network given some noisy and/or bad data 
from the network~\cite{liacco1982role}. 
Estimates are obtained via the (generally non-linear) measurement model:
\begin{equation} \label{eq:state-est}
    \bz = \bh(\bx) + \bepsilon,
\end{equation}
where $\bz \in \mathbb{C}^m$ is a vector of measurements, $\bx \in \mathbb{C}^n$ is a vector of quantities to estimate (typically, voltage phasors), $\bh(\cdot)$ is a vector of functions representing the system physics (i.e., power-flow equations), and $\bepsilon$ is a vector of measurement noise. The state-estimation task is then to estimate $\bx$ given $\bz$ and some knowledge of $\bh(\cdot)$ 
(e.g., its Jacobian matrix).
State estimation has been thoroughly addressed in transmission networks, wherein system  \eqref{eq:state-est} is typically \emph{overdetermined} and \emph{fully observable}: that is, 
(i) the number of measurements $m$ is at least the number of unknowns $n$, and (ii) 
the Jacobian $\bJ \in \comps^{m \times n}$ of $\bh(\cdot)$
is (pseudo) invertible in the sense that 
$
(\bJ^{\sfT} \bJ)^{-1}
$
exists. 
As transmission systems conventionally have redundant measurements that satisfy the observability requirement, 
classical least-squares estimators are applicable and can operate efficiently \cite{Abur2004}.

In contrast, the use of state estimation has historically been limited in distribution networks~\cite{dehghanpour2018survey}. 
Due to limited availability of real-time measurements from Supervisory Control and Data Acquisition (SCADA) systems, 
equation~\eqref{eq:state-est} is typically \emph{underdetermined} ($m < n$), rendering standard least-squares
methods inapplicable. 
Accurate distribution system state estimation was also previously unnecessary since distribution networks only delivered power in one direction towards the customer, requiring minimal distribution system control.
This 
led industry to in practice use only simple heuristics (e.g. based on simple load-allocation rules~\cite{deng2002branch,pereira2004integrated}) to roughly calculate power flow. 

However, due to the increasing adoption of distributed energy resources (DERs) at the edge of the network \cite{driesen2008design},  distribution system state estimation  has become increasingly important \cite{primadianto2017review}. 
There is thus a large focus in the literature on low-observability state estimation techniques.
Many existing methods attempt to improve system observability, 
e.g., by optimizing the placement of additional system sensors \cite{singh2009, bhela2018, YCNN} or by deriving pseudo-measurements from existing sensor data \cite{manitsas2012, wu2013}.
Unfortunately, installation of additional sensors may be expensive or slow, and pseudo-measurements can introduce estimation errors \cite{clements} or be extremely data-intensive to obtain \cite{manitsas2012, YCNN}. 
Other methods seek to perform state estimation using neural networks, without constructing an underlying system model \cite{NNVoltages}.
While such machine learning methods can obtain accurate estimation results, training these methods requires a significant amount of historical data, which may not be available.
As such, there is a need for state estimation methods that can exploit problem structure to perform state estimation at current levels of data availability and observability.

In this paper, we propose a low-observability state estimation algorithm based on \emph{matrix completion} \cite{ExactMatrix}, a tool for estimating missing values in low-rank matrices.
We apply this tool to state estimation 
for a given time step 
by forming a structured data matrix whose rows correspond to measurement locations, and whose columns correspond to measurement types (e.g. voltage or power).
While methods such as \cite{bhela2018} require collecting data over large time windows, our approach enables ``single shot'' state estimation that employs only data from a single time instance.
Our approach is closely related to recent works \cite{gao2016,genes2018} that use matrix completion to estimate lost PMU data over a time series, but while these works estimate missing quantities exclusively at measurement points, we consider the problem of estimating quantities even at non-measurement points where the quantities to be estimated may have never been measured.

The main contributions of our paper are:
\begin{itemize}
    \item A novel distribution system state estimation method based on constrained matrix completion. By augmenting matrix completion with noise-resilient power flow constraints, the proposed method can accurately estimate voltage phasors under low-observability conditions where standard (least-squares) methods cannot.
    \item A flexible framework for employing various types of distribution system measurements into state estimation. Whereas many works (e.g. \cite{klauber, wu2013, bhela2018}) require specific 
    measurements for estimation,
    our approach can accommodate any quantities measured in the field.
    \item An empirical demonstration of the robustness of our method to 
    data availability and measurement loss.
\end{itemize}

The rest of the paper is organized as follows: Section~\ref{sec:mc} introduces the concept of constrained matrix completion. The proposed distribution system state estimation algorithm is presented in Section~\ref{sec:method}. Simulation results for the IEEE 33- and 123-bus test systems are presented in Section~\ref{sec:results}. 
Section~\ref{sec:conclusion} concludes the paper.

\section{Matrix Completion Methods}
\label{sec:mc}

We start by introducing \emph{constrained matrix completion}, a method that is central to our proposed approach.

\subsection{Matrix Completion}
\label{sec:mc-sub}
Given an incomplete matrix that is assumed to be low-rank, the matrix completion problem aims to determine the unknown elements in this matrix.
Formally, let $\M \in \mathbb{R}^{\mcrows \times \mccols}$ be a real-valued 
data matrix, 
$\known \subseteq \{1, \ldots, \mcrows\} \times \{1, \ldots, \mccols\}$ 
describe the
known elements in $\M$, 
and 
$\M_\known \in \mathbb{R}^{\mcrows \times \mccols}$ denote the observation matrix, where $(\M_\known)_{j,k} = \M_{j,k}$ for $(j,k) \in \known$ and $0$ otherwise.
Matrix completion can then be formulated as a rank-minimization problem \cite{ExactMatrix}:
\begin{equation} 
\label{eq:rank-min}
\begin{aligned}
\minimize_{\X \in \reals^{\mcrows \times \mccols}}& \;\; \operatorname{rank}(\X) \\
\subjectto& \; \;\X_{\known}  = \M_{\known},
\end{aligned}
\end{equation}
%
where the decision variable $\X$ estimates $\M$. As the optimization problem~\eqref{eq:rank-min} is NP-hard due to the non-convexity of the rank function, it is common to use a heuristic approach that instead minimizes the \emph{nuclear norm} of the matrix \cite{ExactMatrix}: 
\begin{equation} 
\label{eq:mc}
\begin{aligned}
\minimize_{\X \in \reals^{\mcrows \times \mccols}}& \;\; \|\X\|_* \\
\subjectto& \; \;\X_{\known}  = \M_{\known},
\end{aligned}
\end{equation}
where 
$\|\X\|_*$ sums the singular values of $\X$.
%
Given a sufficient
number of
randomly-sampled
entries in $\M_{\known}$
(depending on the matrix size and rank),
problem~\eqref{eq:mc} often has a unique minimizer $\X$ that 
equals $\M$ \cite{ExactMatrix}.
Additionally,
this problem
can be solved efficiently \cite{TruncNorm,LargeData,AnyBasis}.

Due to the nature of the equality constraint, formulation \eqref{eq:mc} is highly susceptible to noise. 
To alleviate this problem, \cite{candes2010matrix}  proposed an algorithm to handle noisy measurements. 
The algorithm modifies the equality constraint in \eqref{eq:mc} to 
\begin{equation} 
\label{eq:noisy-meas-constr}
\begin{aligned}
\| \X_{\known} - \M_{\known} \|_F \leq \delta,
\end{aligned}
\end{equation}
where $\| \cdot \|_F$ is the Frobenius norm
%
%
and $\delta \geq 0$ is a parameter that can be tuned based on the extent of measurement noise.

\subsection{Constrained Matrix Completion}
\label{sec:constrained-mc}

Now suppose that the values in the matrix 
$\M$ come from some physical system (e.g., a power system).
It is then natural to extend formulation~\eqref{eq:mc}/\eqref{eq:noisy-meas-constr} to incorporate system physics via the following \emph{constrained} optimization problem:
\begin{equation} 
\label{eq:mc_constr}
\begin{aligned}
\minimize_{\X \in \reals^{\mcrows \times \mccols}}& \;\; \|\X\|_* \\
\subjectto& \; \; \|\X_{\known}  - \M_{\known}\|_F \leq \delta, \\
& \; \; \|\bg(\X)\| \leq \beta,
\end{aligned}
\end{equation}
for 
$\delta,\beta \geq 0$, where $\bg(\cdot)$ is a vector of functions 
representing
system physics
(e.g., 
power-flow equations). We note that:
\begin{itemize}
    \item The additional constraint $\|\bg(\bX)\| \leq \beta$ incentivizes low-rank solutions that respect the system physics.
    \item The choice of $\delta$ and $\beta$ is problem dependent. These parameters can be chosen based on the extent of measurement noise, or the objective function can be augmented with terms that try to minimize their values.
    \item If $\bg(\cdot)$ is nonlinear,~\eqref{eq:mc_constr} 
    is non-convex and NP-hard.
\end{itemize}

\section{Low-Observability State Estimation}
\label{sec:method}

We now present our low-observability distribution system state estimation algorithm, which employs the constrained matrix completion 
model
\eqref{eq:mc_constr}.
We describe our power system model, possible formulations of $\M$, and possible 
power system constraints $\bg(\cdot)$
before showing our full formulation.

\subsection{Power System Model}
Let $\buses$ denote the set of buses, where bus 1 is the slack bus and the remaining $\nbus - 1$ buses are $PQ$ buses.
Further, let $\lines \subseteq \buses \times \buses$ denote the set of distribution lines.
We describe the nodal admittance matrix $\Y \in \mathbb{C}^{\nbus \times \nbus}$ in block form as
\[ \Y = \begin{bmatrix*}[l] 
\Y_{\slack \slack} \; \, \scalemath{0.9}{\in \mathbb{C}} \;\;\; & 
\Y_{\slack L} \; \, \scalemath{0.9}{\in \mathbb{C}^{1 \times (\nbus - 1)}} \\ 
\Y_{L \slack} \; \scalemath{0.9}{\in \mathbb{C}^{(\nbus - 1) \times 1}} & 
\Y_{LL} \; \scalemath{0.9}{\in \mathbb{C}^{(\nbus - 1) \times (\nbus - 1)}} 
\end{bmatrix*}.  \]
Let $\vv \in \mathbb{C}^{\nbus}$ and $\sinj \in \mathbb{C}^{\nbus}$ be the vectors of (partially unknown) voltage phasors and net complex power injections, respectively, at each bus.
We denote the slack bus voltage phasor and power injection as $\vv_{\slack}$ and $\sinj_{\slack}$, respectively, and similarly denote the 
vectors of non-slack bus voltages and power injections as $\vvm$ and $\sinjm$. 
Finally, let $\curr \in \mathbb{C}^{\nlines}$ be the vector of (partially unknown) complex currents in each branch, where $\curr_{ft}$ is the current in line $(f,t) \in \lines$.

\subsection{Data Matrix Formulation}

The formulation of the data matrix $\M$ (and thus the optimization variable $\X$) can vary based on the particular attributes of the problem setting, e.g. the kinds of measurements available and problem scale.
We present two possible formulations here, one indexed by branches 
and one indexed by buses.
However, we emphasize that 
\emph{the proposed method is not limited to using these matrix structures}.
The matrix $\M$ can be flexibly structured to accommodate available measurements, as long as these measurements are correlated so that $\M$ is (approximately) low rank.

\subsubsection{Branch Formulation}
\label{sec:branchformu}
$\M$ can be structured
such that each row represents a power system branch and each column represents a quantity relevant to that branch.
This structure allows us to take advantage of both bus- and branch-related measurements.
%
Specifically, for every line $(f, t) \in \lines$, the corresponding row in the matrix $\M \in \mathbb{R}^{\mcrows \times \mccols}$ 
contains:
%
\[
\begin{split}
[\; &\Re(\vv_f), \; \Im(\vv_f), \; |\vv_f|, \; \Re(\sinjp{f}), \; \Im(\sinjp{f}), \; \Re(\vv_t), \; \\
&\Im(\vv_t), \; |\vv_t|, \; \Re(\sinjp{t}), \; \Im(\sinjp{t}), \; \Re(\curr_{ft}), \; \Im(\curr_{ft}) \;], 
\end{split}
\]
where $\mcrows = \nlines$ and we employ $\mccols=12$ quantities per row.

\subsubsection{Bus Formulation}
\label{sec:busformu}
$\M$ can also be structured
such that each row represents a 
bus and each column represents a 
quantity relevant to that bus.
That is, for every
bus $b \in \buses$, the corresponding row in the matrix $\M \in \mathbb{R}^{\mcrows \times \mccols}$ 
contains:
\[
\begin{split}
[\; &\Re(\vv_b), \; \Im(\vv_b), \; |\vv_b|, \; \Re(\sinjp{b}), \; \Im(\sinjp{b})], 
\end{split}
\]
where $n_1 = |\buses|$ and we employ $n_2=5$ quantities per row.
While this structure only employs bus-related measurements,
its advantage is that it yields small matrices that can be used for efficient estimation on large-scale problems.

\subsubsection{On the Low-rank Assumption} 
\label{sec:lowrankassum}

We note that the complex power system quantities in both the bus and branch formulations are approximately linearly correlated, as implied by the fact that the power system equations employing them can be expressed in (approximate) linear form (see Section~\ref{sec:power-sys-constr}).
In other words, these data matrix formulations should be approximately low-rank.

\begin{figure*}[ht!]
    \centering
    \subfloat[IEEE 33-bus feeder (branch formulation). The largest singular value (out of 12) comprises over 98\% of the singular value sum.]{\includegraphics[width=0.47\textwidth]{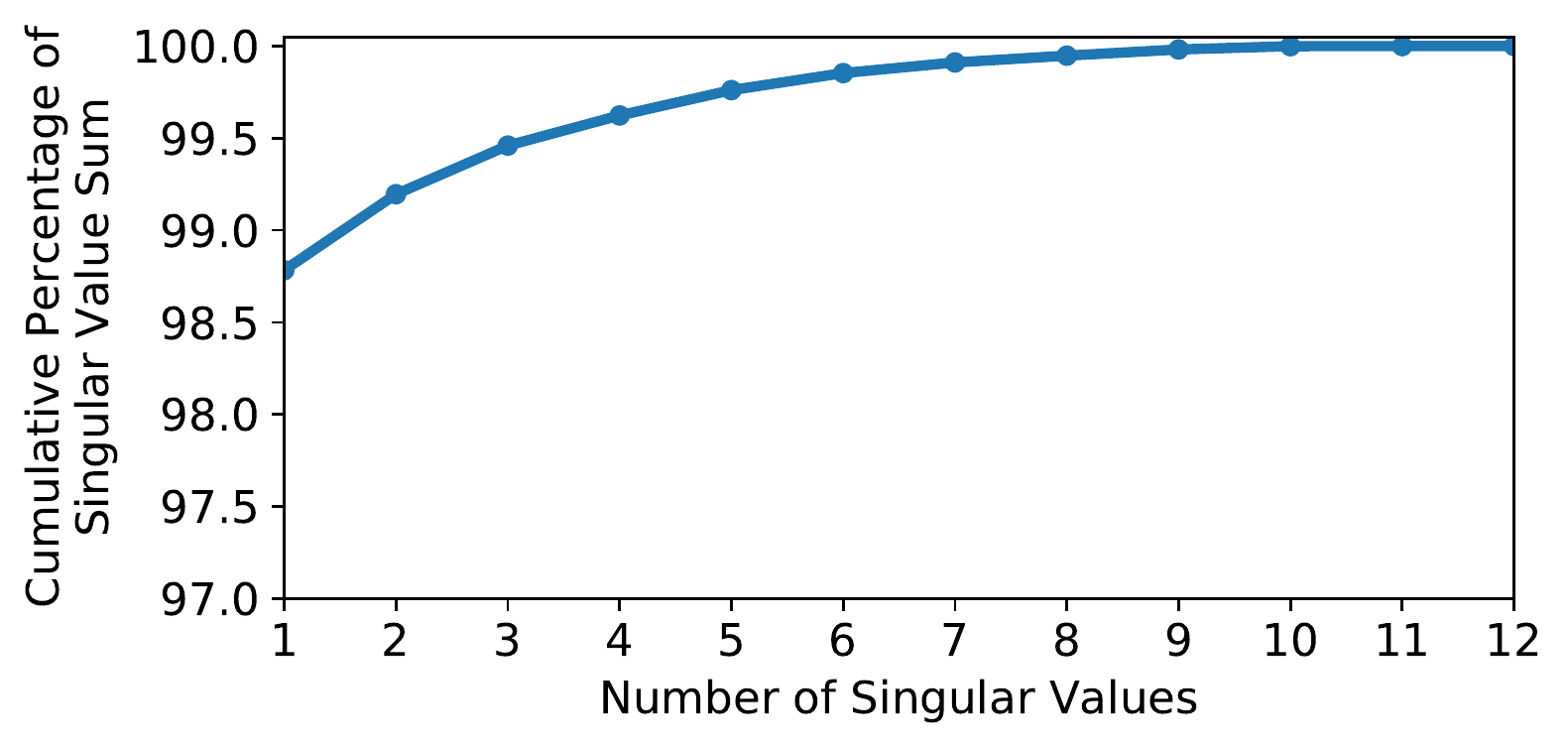}}\hfill
	\subfloat[IEEE 123-bus feeder (multi-phase bus formulation). The first 3 (out of 5) largest singular values comprise 95\% of the singular value sum. ]{\includegraphics[width=0.46\textwidth]{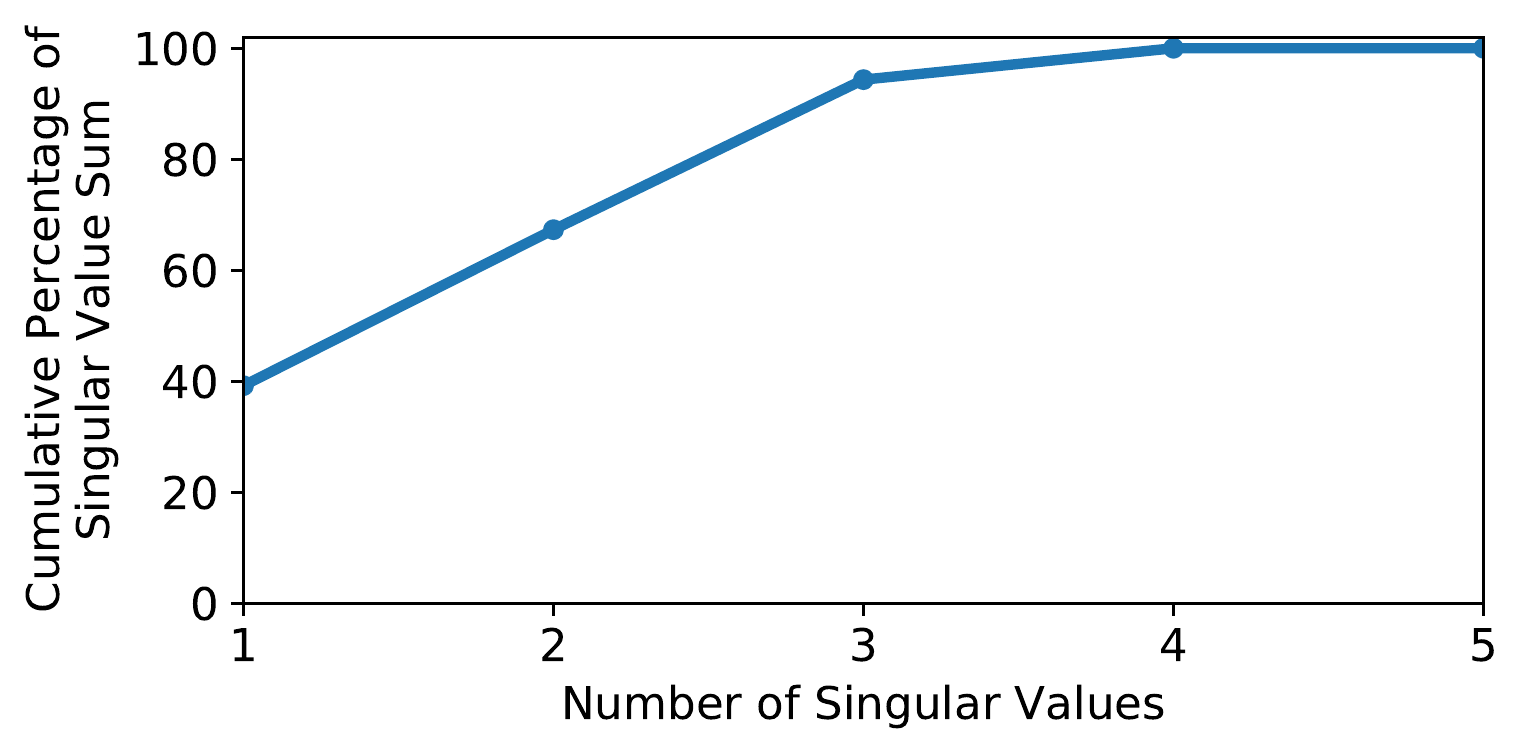}}
    \caption{Cumulative percentage of matrix singular value sum for two test cases. Both matrices are approximately low-rank.}
    \label{fig:sing-vals}
\end{figure*}
We observe empirically that 
this low-rank 
assumption holds in practice.
Fig.~\ref{fig:sing-vals} shows the cumulative percentage 
distribution 
of singular values 
for the IEEE 33-bus feeder (using a branch formulation) and the IEEE 123-bus feeder (using a bus formulation).
In both cases, we see that a few singular values comprise much of the singular value sum, implying that these matrices are (approximately) low-rank.

\subsection{Physical Power Flow Constraints}
\label{sec:power-sys-constr}

As described in Section~\ref{sec:constrained-mc}, we augment 
matrix completion 
with power system constraints to 
encourage 
physically meaningful solutions.
These constraints are linear to ensure that problem~\eqref{eq:mc_constr} is convex.
We describe the constraints we use below, but note that constraints can be added, removed, or modified depending on the types of measurements in $\M$.

\subsubsection{Duplication Constraint}
\label{sec:dup-constr}
Depending on the formulation, some quantities may appear in more than one location in $\M$.
For example, in our branch formulation, 
quantities related to a given bus 
appear in multiple rows if 
the
bus is in multiple branches.
We thus constrain equivalent quantities in the matrix to be equal.
Formally, let $\Lambda$ contain 
all pairs of indices of duplicated quantities in $\M$.
We require that
\begin{equation}
\label{eq:dups}
    \X_{\lambda_1} = \X_{\lambda_2},\; \forall\; (\lambda_1, \lambda_2) \in \Lambda.
\end{equation}

\subsubsection{Ohm's Law Constraint}
\label{sec:ohms-constr}
When $\M$ contains both bus- and branch-related quantities (as in the branch formulation), we can apply Ohm's Law, defined as 
\begin{equation}
\label{eq:ohms}
(\vv_f - \vv_t) \by_{ft} = \curr_{ft}, \; \forall\; (f,t) \in \lines,
\end{equation}
where $\by_{ft}$ is the line admittance.
However, 
using an
exact equality constraint may cause the matrix completion problem to become infeasible,
e.g.
due to measurement noise.
We thus employ a noise-resilient version of Ohm's Law, i.e.,
\begin{equation}
\label{eq:ohms-noise}
\begin{bmatrix}
-\ohmsparam_{r, ft} \\[2pt] -\ohmsparam_{c, ft}
\end{bmatrix} 
\leq 
\begin{bmatrix}
\Re\left((\vv_f - \vv_t) \by_{ft} - \curr_{ft} \right) \\[2pt]
\Im\left((\vv_f - \vv_t) \by_{ft} - \curr_{ft}\right)
\end{bmatrix}
\leq 
\begin{bmatrix}
\ohmsparam_{r, ft} \\[2pt] \ohmsparam_{c, ft}
\end{bmatrix},
\end{equation}
where $\ohmsparam_{r, ft},\; \ohmsparam_{c, ft} \in \mathbb{R}_+$ are respective error tolerances for the real and complex parts of Ohm's Law on line $(f,t) \in \lines$.

\subsubsection{Linearized Power Flow Constraints}
\label{sec:pf-constr}

As the exact AC power flow equations are non-linear, we employ Cartesian linearizations of these equations.
For non-slack voltages and power injections, we employ approximations of the form
\begin{subequations} \label{eq:lin_app}
\begin{align}
\vvm  \; &\approx \bA \begin{bmatrix} \Re(\sinjm) \\ \Im(\sinjm) \end{bmatrix} + \: \zloadv,  \label{eq:vlin} \\
|\vvm| &\approx \bC \begin{bmatrix} \Re(\sinjm) \\ \Im(\sinjm) \end{bmatrix}  +|\zloadv|. \label{eq:vmaglin}
\end{align}
\end{subequations}
For example, using the method proposed in \cite{Bernstein_Dall’Anese_2017}, we can let $\zloadv = -\vv_{\slack} \Y_{LL}^{-1} \Y_{L\slack} \in \mathbb{C}^{\nbus-1}$ be the vector of non-slack zero-load voltages and $\bA,\; \bC \in \mathbb{C}^{(\nbus-1) \times 2(\nbus-1)}$ be defined for some non-slack voltage estimates $\hat{\vv}_{-\slack}$ as
\begin{subequations}
\begin{align}
\bA &= \begin{bmatrix} \Y_{LL}^{-1} \diag(\overline{\hat{\vv}}_{-\slack})^{-1} & -j\Y_{LL}^{-1} \diag(\overline{\hat{\vv}}_{-\slack})^{-1} \end{bmatrix} , \\
\bC &= \diag(|\hat{\vv}_{-\slack}|)^{-1} \Re \left(\diag(|\overline{\hat{\vv}}_{-\slack}|) \bA \right).
\end{align}
\end{subequations}
For our case, 
we let $\hat{\vv}_{-\slack} = \zloadv.$ 
We note, however, that other methods to obtain the linear approximations \eqref{eq:lin_app} 
(e.g., data-driven regression methods~\cite{Liu2017})
can also be leveraged.

To relate voltages with the power injection at the slack bus, we employ the exact power flow equation 
\begin{equation}
\label{eq:slackpow}
\sinjp{\slack} = \vv_{\slack} (\overline{\Y}_{\slack \slack} \overline{\vv}_{\slack} + \overline{\Y}_{\slack L} \overline{\vv}_{-\slack}).
\end{equation}
This equation is linear in voltages since
$\vv_{\slack}$
is known.

As in Section~\ref{sec:ohms-constr}, we relax these constraints into noise-resilient versions as 

\begin{subequations}
\begin{footnotesize}
\noindent
\begin{align} 
\begin{bmatrix} -\vparam_{r} \\[10pt] -\vparam_{c} \end{bmatrix} \; &\leq \,\, 
\begin{bmatrix} 
\Re\Big(\vvm  \, - \Big(\bA \begin{bmatrix} \Re(\sinjm) \\ \Im(\sinjm) \end{bmatrix} + \: \zloadv \, \Big) \Big) \\[10pt]
\Im\Big(\vvm  \, - \Big(\bA \begin{bmatrix} \Re(\sinjm) \\ \Im(\sinjm) \end{bmatrix} + \: \zloadv \, \Big) \Big)
\end{bmatrix} 
\leq
\begin{bmatrix} \vparam_{r} \\[10pt] \vparam_{c} \end{bmatrix} ,
 \label{eq:vlin-noise} \\[10pt]
-\vmagpar \; &\leq |\vvm| - \left(\bC \begin{bmatrix} \Re(\sinjm) \\ \Im(\sinjm) \end{bmatrix}  +|\zloadv|\right) \; \leq \vmagpar, \label{eq:vmaglin-noise} \\[10pt]
\begin{bmatrix} -\alpha_{r} \\[2pt] -\alpha_{c} \end{bmatrix}
&\leq
\begin{bmatrix}
\Re\big(\sinjp{\slack} - \left( \vv_{\slack} (\overline{\Y}_{\slack \slack} \overline{\vv}_{\slack} + \overline{\Y}_{\slack L} \overline{\vv}_{-\slack}) \right)  \big) \\[2pt]
\Im\left(\sinjp{\slack} - \left( \vv_{\slack} (\overline{\Y}_{\slack \slack} \overline{\vv}_{\slack} + \overline{\Y}_{\slack L} \overline{\vv}_{-\slack}) \right)  \right)
\end{bmatrix}
\leq
\begin{bmatrix} \alpha_{r} \\[2pt] \alpha_{c} \end{bmatrix},
\label{eq:slackpow-noise}
\end{align}
\end{footnotesize}
\end{subequations}

\noindent where $\vparam_{r}, \vparam_{c}, \vmagpar \in \mathbb{R}_+^{\nbus-1}$, $\alpha_{r}, \alpha_{c} \in \mathbb{R}_+$ are error tolerances, and inequalities are evaluated elementwise.

\subsection{Full Problem Formulation}
\label{fullproblemformu}
Given these power flow constraints, we collect our error tolerances into the set $\mathcal{T} = \{\ohmsparam_{r}, \ohmsparam_{c}, \vparam_{r}, \vparam_{c}, \vmagpar, \alpha_{r}, \alpha_{c}\}$ and form our constrained matrix completion problem~\eqref{eq:mc_constr} as
\begin{subequations}
\label{eq:mc-pf}
\begin{align} 
\minimize_{\X \in \mathbb{R}^{\mcrows \times \mccols},\; \mathcal{T}}&\;\;
\| \X \|_{*} + \sum_{\bt \in \mathcal{T}} w_{\bt} \| \bt \| \\
\subjectto& \;\; ~\lVert \X_{\known} - \M_{\known} \rVert_F \leq \delta,\\
\qquad& \;\; \eqref{eq:dups}, 
\eqref{eq:ohms-noise}, \eqref{eq:vlin-noise}, \eqref{eq:vmaglin-noise}, \eqref{eq:slackpow-noise}, \\
\qquad& \;\; \bt \geq 0, \;\; \forall \bt \in \mathcal{T},
\end{align}
\end{subequations}
where in this case we add each constraint tolerance $\bt \in \mathcal{T}$ to the objective with an associated 
weight $w_{\bt}$. 
Each weight is chosen to reflect the relative importance of its constraint.
For a branch-formulated $\M$, the above formulation can be used as-is.
For a bus-formulated $\M$, equations~\eqref{eq:dups} and~\eqref{eq:ohms-noise} are removed from the constraints (and their associated parameters $\ohmsparam_{r}, \ohmsparam_{c}$ are removed from $\mathcal{T}$) since $\M$ does not contain duplicated quantities or branch current measurements.

Formulation~\eqref{eq:mc-pf} allows the matrix completion optimization 
to explicitly trade off between the low-rank assumption and fidelity to the power flow constraints, without requiring tuning of each entry of each constraint tolerance vector.
We further observe that, since the objective is convex and all constraints are linear in the entries of $\X$, this formulation is a convex optimization problem and can be solved efficiently.

\subsection{Extension to the Multi-phase Setting}
While for brevity we formally present only a single-phase, balanced formulation, our approach can easily be extended to the 
general multi-phase setting.
Specifically, $\M$ can be structured to include phase-wise quantities, and the constraints presented in Section~\ref{sec:power-sys-constr} can be replaced with multi-phase versions (e.g. see~\cite{Bernstein_Dall’Anese_2017, bernstein2017}).
We empirically illustrate the application of our method to a multi-phase $123$-bus test case in Section \ref{sec:123_bus}.

\section{Simulation and Results}
\label{sec:results}
We demonstrate the performance of our matrix completion method on the IEEE 33-bus and 123-bus test cases.
We employ the branch-formulation of our algorithm on the 33-bus feeder,
and show that it performs well in low-observability settings (where traditional state estimation techniques cannot operate) as well as in full-observability settings.
We also evaluate our method on the multi-phase 123-bus feeder using a bus formulation matrix,
demonstrating that our method scales robustly to larger systems.

\subsection{33-Bus System}

We test our 
branch-formulated 
algorithm on a modified version of the IEEE 33-bus test case with solar panels added at buses 16, 23, and 31.
On this system, voltage magnitudes range from approximately $0.99$-$1.02$ p.u., and all angles (relative to the substation voltage angle) are close to 0.
We assume that voltage phasors are known \emph{only} at the slack bus, and must be estimated elsewhere. 
Non-voltage phasor quantities (i.e.\ voltage magnitude, power injections, and current flows) are assumed to be known exactly at the feeder head, and are ``potentially known'' at other buses.

\subsubsection{Randomly-sampled data}

\begin{figure}
    \centering
    \includegraphics[width=0.46\textwidth]{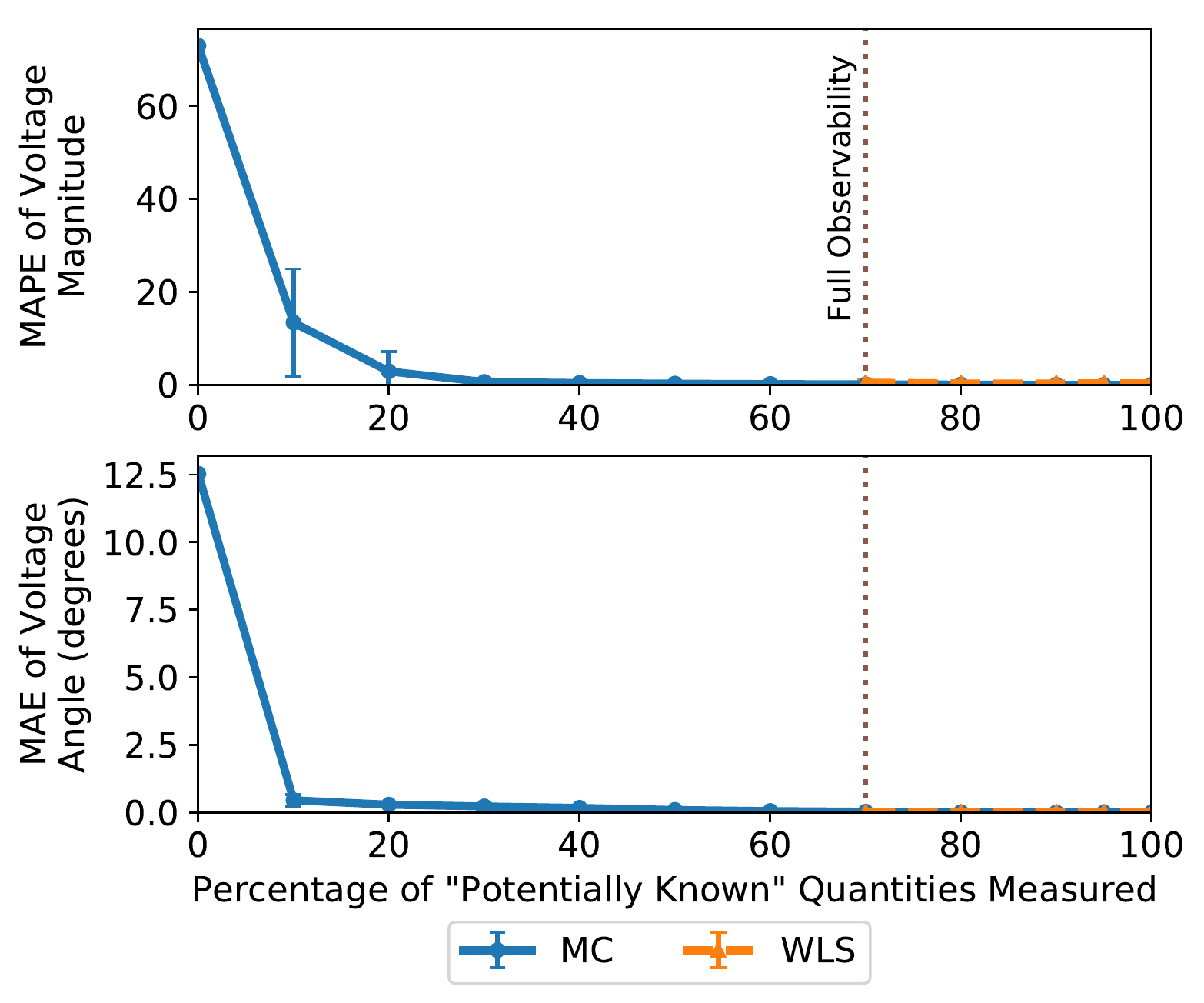}
     \caption{Performance on the 33-bus test case with 1\% noise and random sampling of data. (Each point represents 50 runs.) We achieve less than $1\%$ voltage magnitude MAPE and less than $0.23$ degrees voltage angle MAE when $30\%$ or more of ``potentially known'' quantities are measured.}
    \label{fig:33bus-randsamp-noisy-summary}
\end{figure}

In one set of experiments, we model data unavailability among ``potentially known'' quantities via random sampling.
That is, we randomly choose sets of buses (ranging between $0$-$100\%$ of buses) at which all ``potentially known'' quantities are known, and set these quantities to be unknown at all other buses.
Results under $1\%$ Gaussian measurement noise are shown in Fig.~\ref{fig:33bus-randsamp-noisy-summary}.
In both cases, the mean absolute percent error (MAPE) of our voltage magnitude estimates drops to below $3\%$ and the mean absolute error (MAE) of our voltage angle estimates drops to below $0.29$ degrees when $20\%$ or more measurements are known.
(We report MAEs rather than MAPEs for voltage angles since all angles are close to 0.)
When $30\%$ or more of measurements are known, voltage magnitude MAPEs are much less than $1\%$, and voltage angle MAEs are less than $0.23$ degrees.
This result calibrates well with guarantees for matrix completion performance under random data removal \cite{candes2010matrix}.
In contrast, traditional full-observability state estimation techniques cannot 
operate on this system until $70\%$ of ``potentially known'' quantities are measured.
At full observability, our method's estimates are competitive with those of a state-of-the-art
weighted least-squares (WLS)
state estimation algorithm,
with voltage magnitude MAPEs and angle MAEs within $0.2\%$ and $0.03$ degrees, respectively, for  our  matrix  completion  algorithm  and  within  $0.5\%$  and $0.008$  degrees,  respectively,  for  WLS. 

\begin{figure*}[t!]
    \centering
    \includegraphics[width=\textwidth]{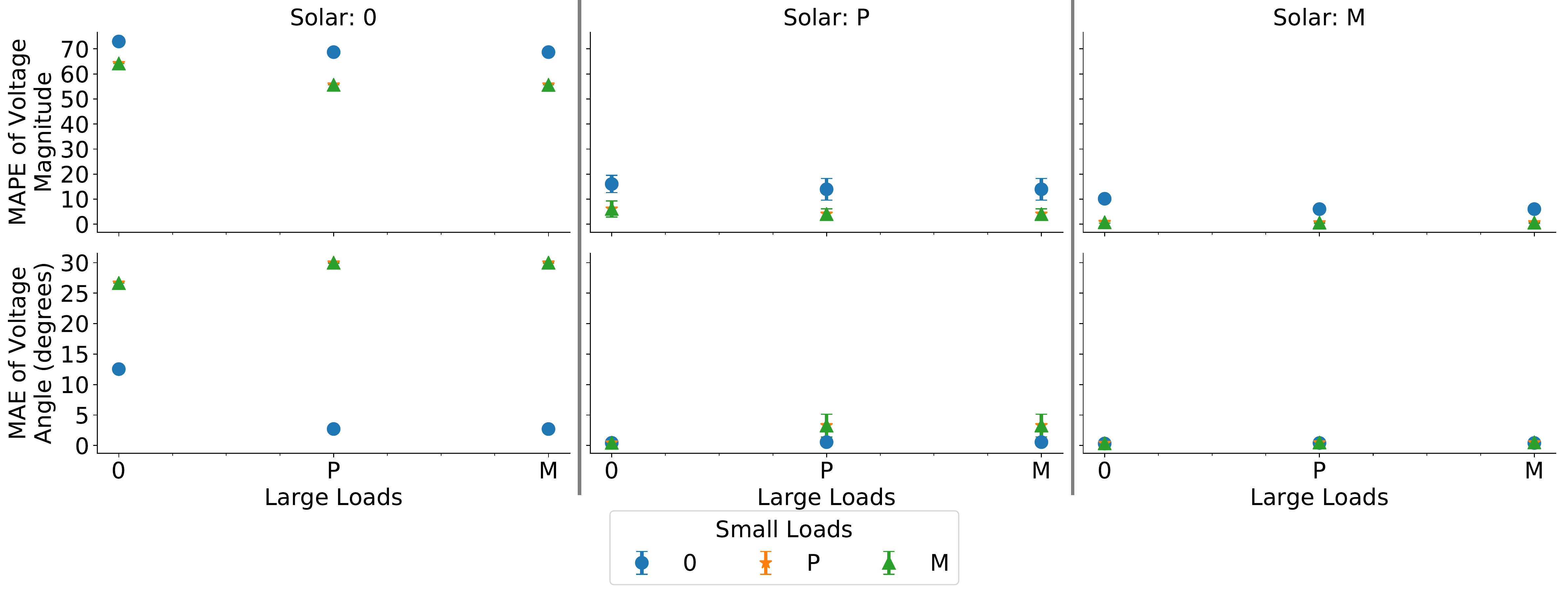}
    \caption{Performance on the 33-bus test case with $1\%$ noise over different data availability for solar, large load, and small load buses (0 = not measured, P = some pseudomeasurements, M = some measurements; 
    each point represents 50 runs.) Our 
    method achieves good estimation error in many data availability scenarios, even though all scenarios exhibit low observability.}
    \label{fig:33bus-datadriven-noisy-summary}
\end{figure*}
\begin{figure}[ht!]
    \centering
    \includegraphics[width=0.46\textwidth]{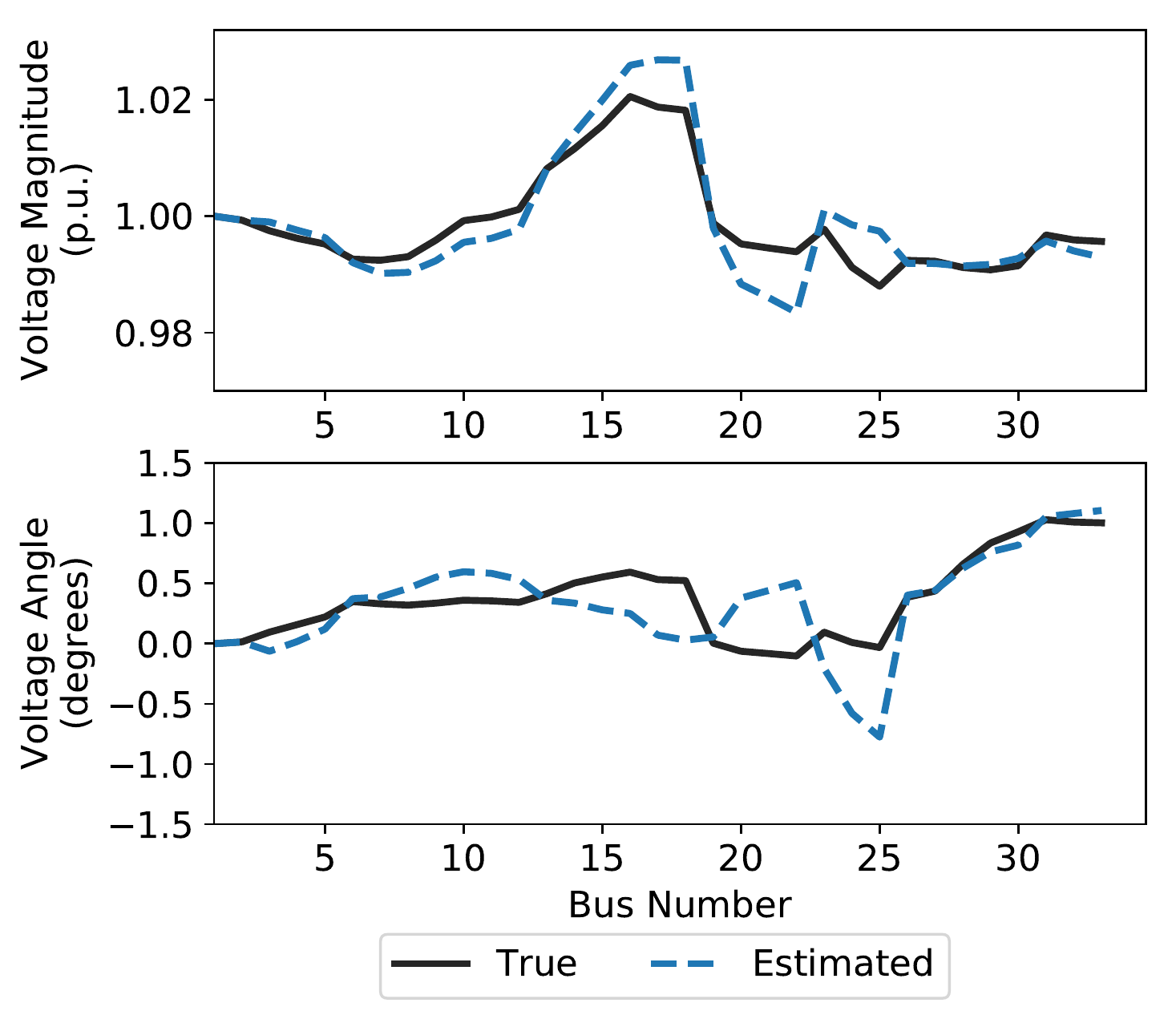}
    \caption{Representative voltage phasor estimates for the 33-bus test case with $1\%$ noise in the low-observability scenario 
    with solar, large loads, and small loads partially measured. 
    }
    \label{fig:33bus-datadriven-noisy-rep}
\end{figure}

\subsubsection{Data-driven assumptions}

In practice,
data unavailability is not uniformly random, but instead systematic and correlated.
For instance, a utility may only have certain types of sensors 
at certain types of buses. 
We thus
run a second set of experiments where we classify the buses into four categories: slack bus (1 bus), solar generators (3 buses), large loads (6 buses), and small loads (17 buses). 
As before, all quantities are known exactly at the slack bus.
At solar PV generators, 
real power injections, reactive power injections, and voltage magnitudes are potentially known.
At loads (large or small), real power injections are potentially known.
Since actual sensor availability may vary between utilities,
we model different scenarios in which these groups of non-slack buses have real, pseudo-, or no measurements.
In the best case (when all three groups of buses have some measurements), $23\%$ of ``potentially known'' quantities are measured.
Thus, \emph{all our data-driven scenarios are at low observability}, and traditional full-observability state estimation methods cannot be used.

The performance of our method on of these scenarios is shown in Fig.~\ref{fig:33bus-datadriven-noisy-summary}, assuming measurements have $1\%$ Gaussian sensor noise
and pseudomeasurements have $10\%$ Gaussian error.
Our algorithm 
achieves less than $1\%$ MAPE in its magnitude estimates and less than $0.5$ degrees MAE in its angle estimates 
when accurate measurements are available for solar generators and either measurements or pseudomeasurements are available for loads. 
Results for a representative run in this scenario are shown in Fig.~\ref{fig:33bus-datadriven-noisy-rep}.
More generally, 
our estimates (averaged across all runs) have at most $10.2 \pm 0.2$\% voltage magnitude MAPE and $0.50 \pm 0.26$ degrees voltage angle MAE in \emph{any} scenario where solar generators are measured, and at most $16.0 \pm 3.4$\% magnitude MAPE and $3.3 \pm 1.9$ degrees angle MAE in any scenario where solar generators have pseudomeasurements (where angle errors are high in this latter case due to solar generator measurement noise).
If solar generator measurements are unknown, magnitude MAPEs range from $57$-$73\%$, and angle MAEs range from $2.7$-$30$ degrees.

\begin{figure}[t!]
    \centering
    \includegraphics[width=0.46\textwidth]{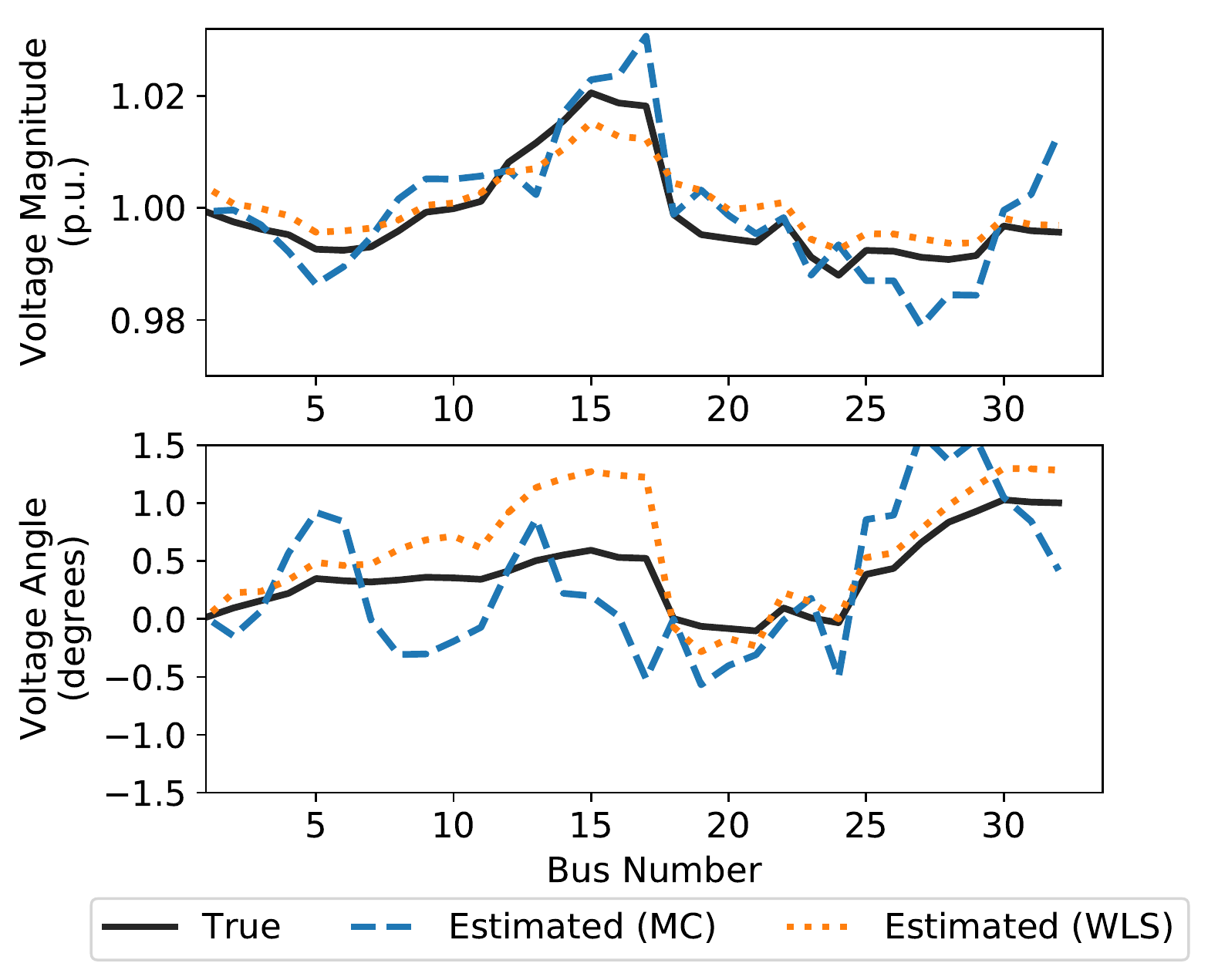}
    \caption{Representative 
    estimates for the 33-bus test case with $1\%$ noise and data-driven sampling at full observability.
    }
    \label{fig:33bus-datadriven-fullrank-noisy-rep}
\end{figure}

A potential alternative to using low-observability state estimation techniques is to enable full-observability 
techniques
by deploying additional 
sensors.
To model this alternative, we randomly add AMI sensors (which collect coarse-granularity load data, modeled as pseudomeasurements) and ``magnitude sensors'' (which collect voltage and current magnitudes) to our system until full observability is achieved.
We then compare our method to WLS
on this augmented system.
Voltage phasor estimates for a representative run are shown in Fig.~\ref{fig:33bus-datadriven-fullrank-noisy-rep}.
In this case, both our matrix completion algorithm and WLS quite accurately estimate voltages, 
with voltage magnitude MAPEs and angle MAEs within $0.6\%$ and $0.53$ degrees, respectively, for our matrix completion algorithm and within $0.8\%$ and $0.67$ degrees, respectively, for WLS.
However, 
achieving full observability required
adding 28-63 sensors to the system (depending on the baseline data availability scenario), which represents a potentially high cost to the distribution utility.

Overall, our results on the 33-bus system demonstrate that our matrix completion algorithm can provide accurate state estimation performance in the low-observability case (where traditional state estimation techniques cannot operate), as well as in the full-observability case.

\subsection{123-Bus Feeder} 
\label{sec:123_bus}


We next demonstrate that our matrix completion method effectively scales to larger systems via experiments on the IEEE 123-bus feeder.
The 123-bus feeder is a multi-phase unbalanced radial distribution system, in which buses are single-, double-, or three-phase (with 263 phases in total.) 
On this system, voltage magnitudes range from $0.95$-$1$ p.u., and voltage angles are around $0$ or $\pm 120$  degrees.
We employ a bus-formulation matrix for this test case,
where each matrix row represents a phase at one bus.
To validate our approach, we use two hours of system voltage and power injection data.
This data was simulated at one-minute resolution using power flow analysis with diversified load and solar profiles created for each bus.

As in the previous section, we assume that voltage phasors are known \emph{only} at the slack bus (which here has three phases) and must be estimated elsewhere.
All other measurements (i.e.\ voltage magnitudes and power injections) are known at the slack bus and are ``potentially known'' for non-slack system phases.
\begin{figure}[t!]
    \centering
    \includegraphics[width=0.46\textwidth]{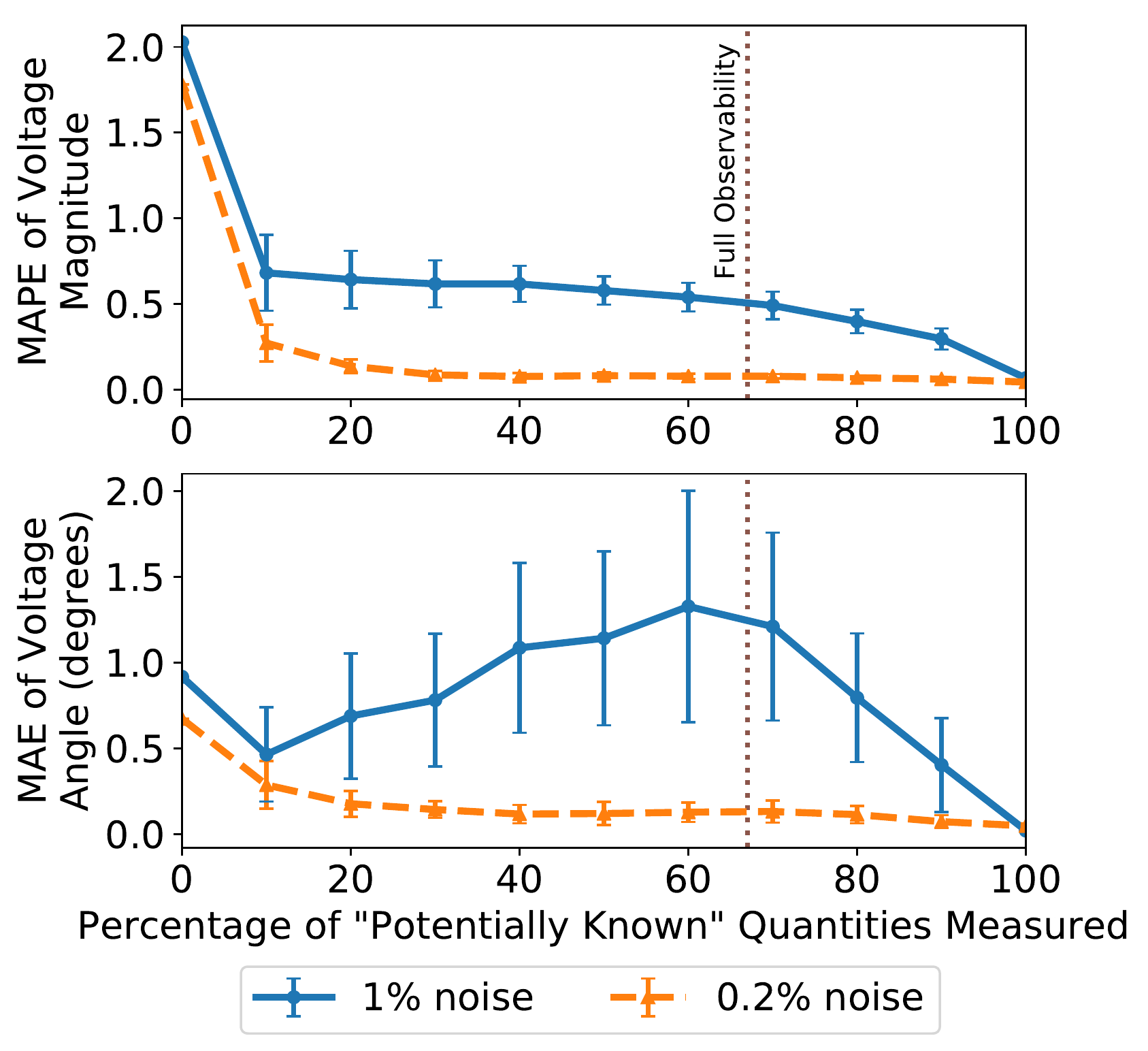}
    \caption{Performance on the 123-bus test case for one time step.
    (Each point represents 50 runs.)
    We achieve less than 1\% voltage magnitude MAPE when 10\% or more of ``potentially known'' quantities are measured.
    The voltage angle MAE is always below $1.5$ degrees at $1\%$ measurement noise, and below $1$ degree at $0.2\%$ measurement noise.
    }
    \label{MAPE_volmagangcomp}
\end{figure}

\begin{figure}[t!]
    \centering
    \includegraphics[width=0.46\textwidth]{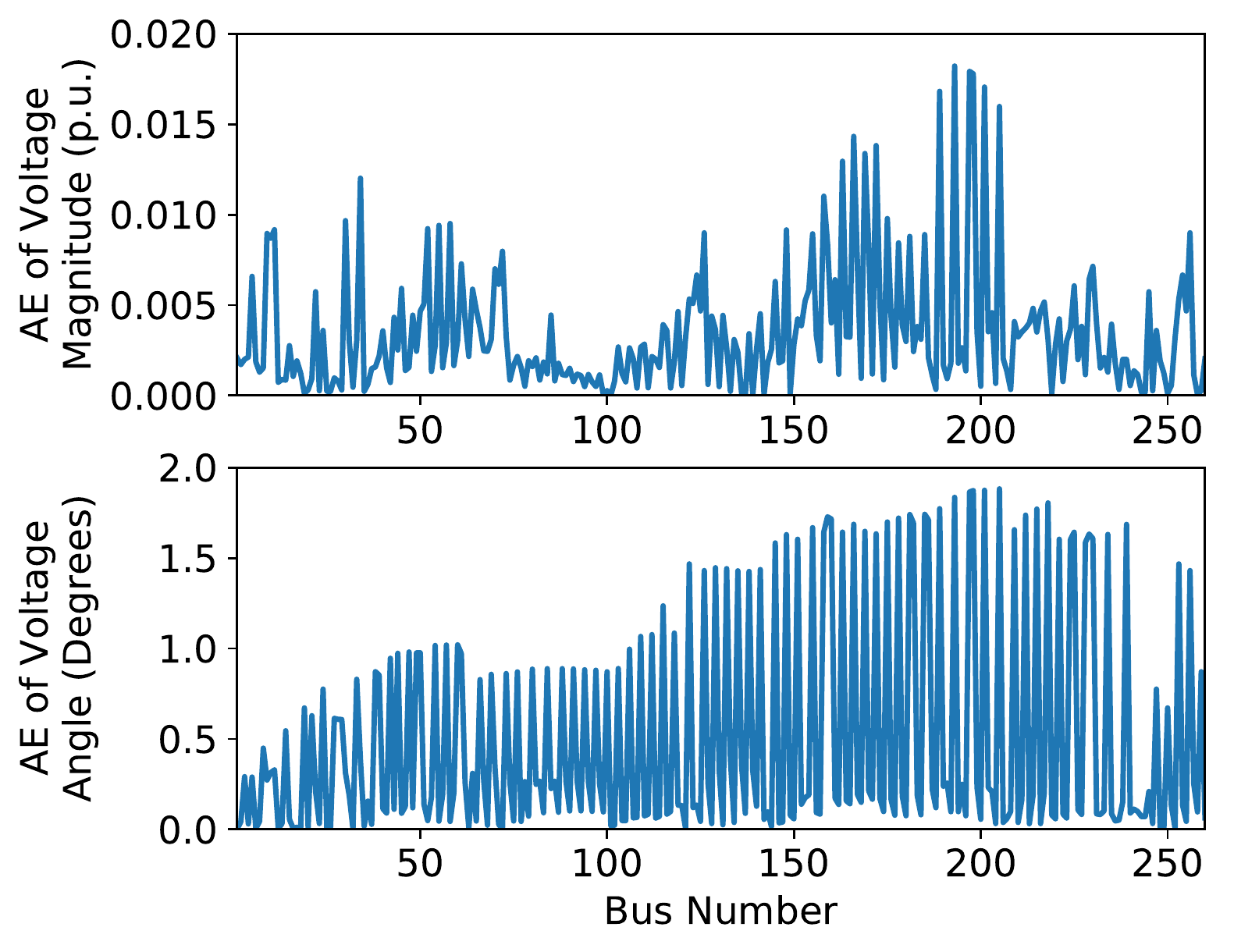}
    \caption{Performance on the 123-bus test case for a representative run 
    with $50\%$ data availability
    and $1\%$ noise. 
    }
    \label{fig:123bus-50pct-abserr}
\end{figure}

We first employ our algorithm at one point in the time series during which solar injections are nonzero.
In our experiments, we vary the percentage of ``potentially known'' quantities that are measured, and note that this system exhibits low-observability if less than two-thirds of these quantities are measured.
Results for the cases of $0.2\%$ and $1\%$ measurement noise are shown in Fig.~\ref{MAPE_volmagangcomp}.
We also show representative results for one run under $1\%$ measurement noise and $50\%$ measurement availability (which is in the low-observability realm) in Fig.~\ref{fig:123bus-50pct-abserr}.

These results show that our algorithm estimates voltage phasors with relatively high accuracy across all levels of data availability.
The MAPE of our voltage magnitude estimates is less than $2.6\%$ even when no ``potentially known'' quantities are measured, and falls to less than $1\%$ once $10\%$ or more of these quantities are available.
Unsurprisingly, our voltage magnitude estimates are better when measurement noise is lower.
We do note, however, that the MAPE is not equal to zero even when $100\%$ of measurements are available, since we use approximate linear power flow equations as constraints in our formulation~\eqref{eq:mc-pf}.

For voltage angle estimates, we again report MAEs rather than MAPEs since the angles at some phases are close to zero.
For $1\%$ measurement noise, we see that the MAE is at most $1.75$ degrees across all data availability levels, which is small given that most voltage angles are around $\pm 120$ degrees. 
However, the accuracy of our angle estimates does not decrease monotonically as more measurements are added.
This is because our algorithm directly estimates the real and imaginary parts of voltage phasors; while we estimate these quantities accurately, their errors are not correlated (especially under large measurement noise), leading to inaccuracies in their implied angle estimates.
At $0.2\%$ measurement noise, the MAE for voltage angle decreases as more measurements are available, dropping below 1 degree even when no measurements are available.

\begin{figure}[t!]
		\centering
        \includegraphics[width=0.46\textwidth]{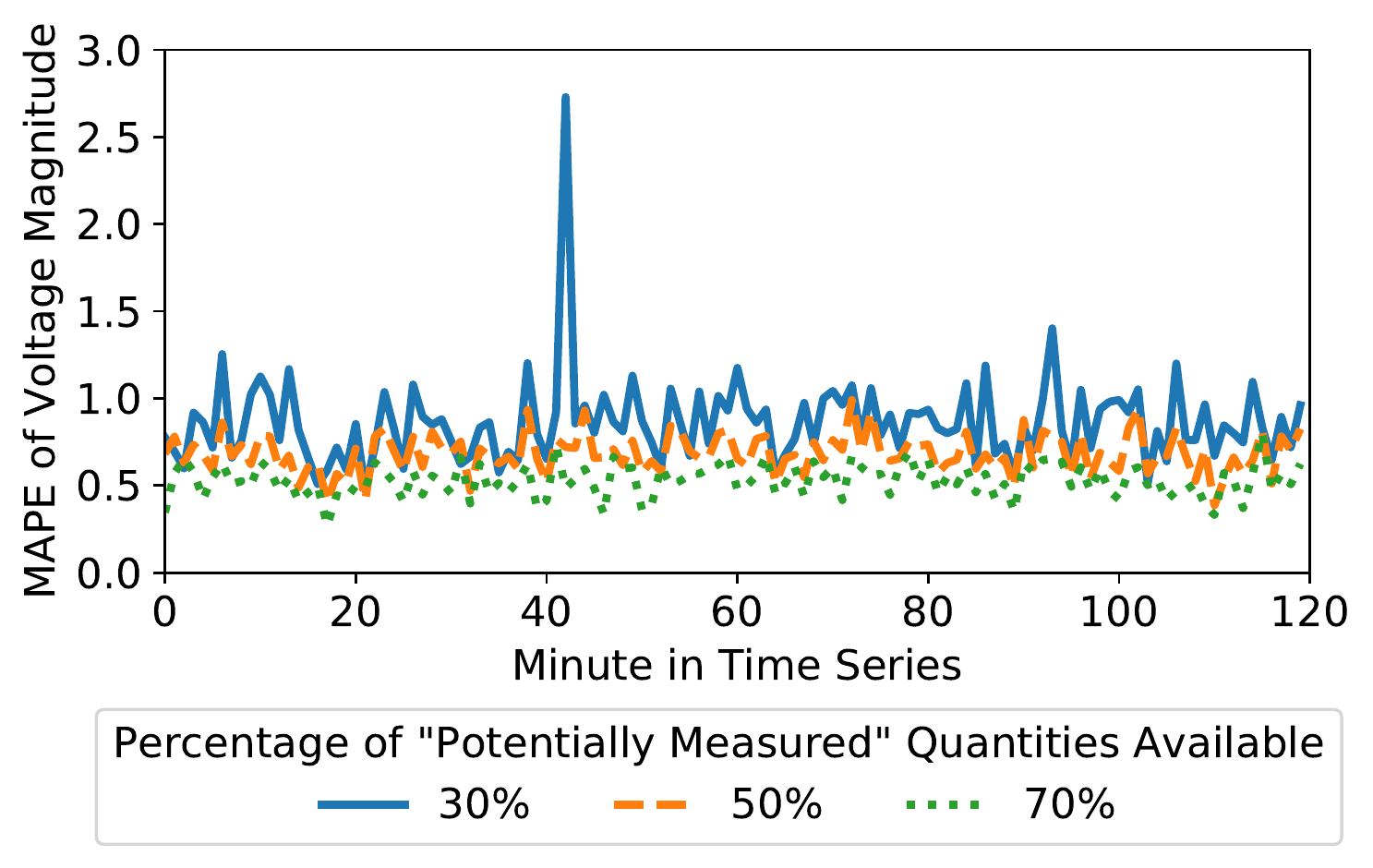}
        \caption{Time-series performance for the 123-bus test case with $1\%$ noise.
        Voltage magnitude MAPEs are at most $2.7\%$ with $30\%$ of data, and at most $1\%$ with $50\%$ or $70\%$ of data.
		}
       \label{MAPE_volmag_twohour}
\end{figure}

To demonstrate the scalability of our method, we next implement our matrix completion algorithm on the entire two-hour time series. 
We model data availability by placing sensors at fixed sets of buses comprising $30\%$, $50\%$, or $70\%$ of all buses (with larger sets of buses inclusive of smaller sets).
The MAPEs of our voltage magnitude estimates under 1\% measurement noise are shown in Fig.~\ref{MAPE_volmag_twohour}.
We find that under $30$\% data availability, the MAPEs of these estimates are within $1\%$ for $80\%$ of time steps, and within $1.5\%$ for $99\%$ of time steps, with a maximum MAPE of $2.7\%$. 
Under $50\%$ and $70\%$ data availability, the MAPEs of our voltage magnitude estimates are within $1\%$ across all points in the time series. 
The MAEs of our voltage angle estimates at each time step (not pictured) are similar to those in Fig.~\ref{MAPE_volmagangcomp}.

\begin{figure}[t!]
		\centering
        \includegraphics[width=0.46\textwidth]{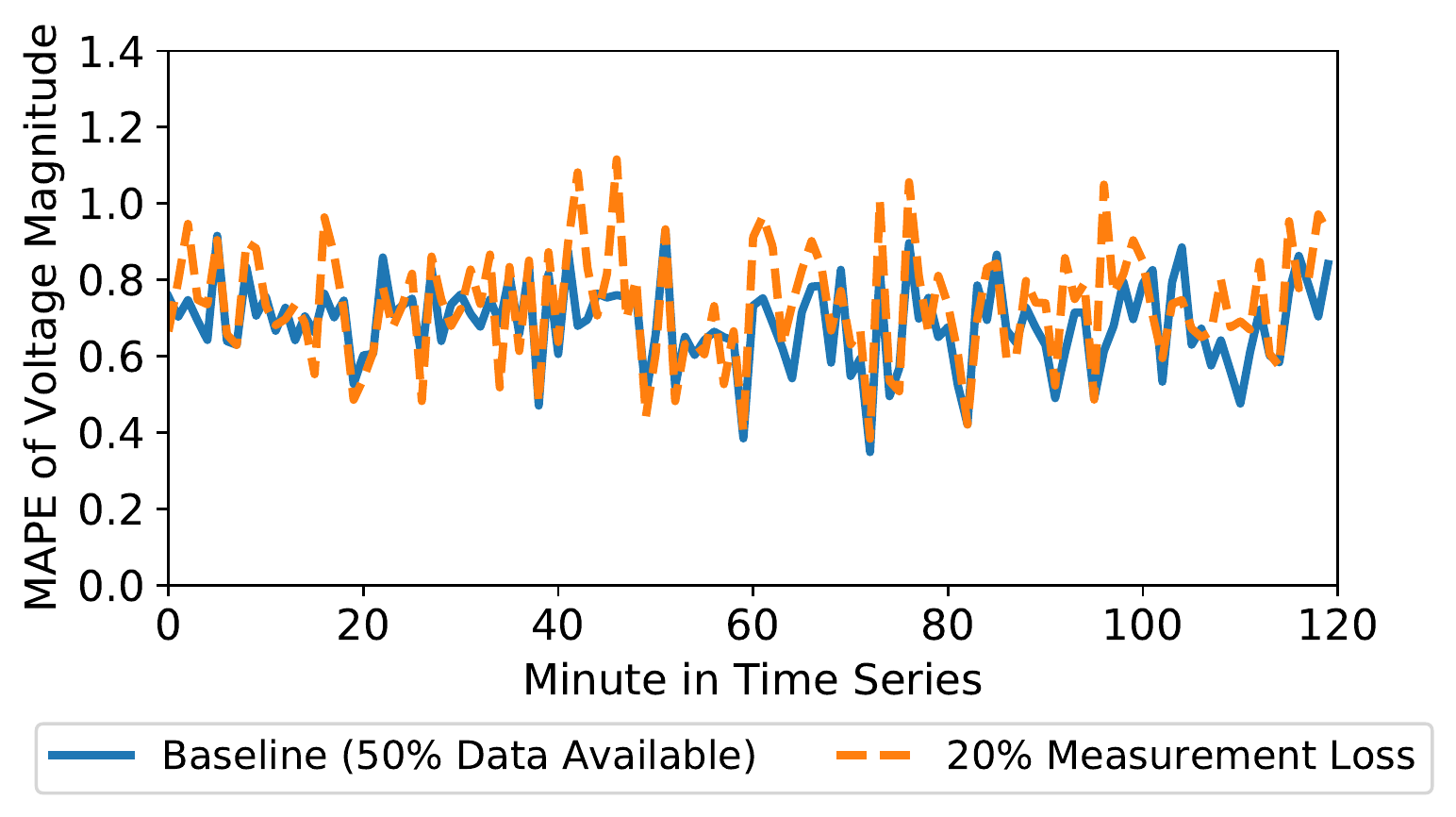}

        \caption{Time-series performance for the 123-bus test case with $1\%$ noise when measurements are randomly lost. 
        Our magnitude estimates are relatively robust to data loss.
        }
        \label{MAPE_loss_twohour}
\end{figure}

Finally, we demonstrate the robustness of our algorithm to dynamic measurement loss, which commonly occurs on the distribution system due to failures such as packet drops and equipment malfunctions.
We model the baseline data availability as being $50\%$ (low-observability), and randomly remove $20\%$ of the available measurements at each point in the two-hour time series.
The MAPEs of our voltage magnitude estimates under $1\%$ measurement noise are shown in Fig.~\ref{MAPE_loss_twohour}.
We find that in most cases, our estimates under data loss are similar to the estimates without data loss, with larger deviations in some cases when critical measurements are lost.
In all cases, the MAPEs of our voltage magnitude estimates are less than $1.2\%$, demonstrating the robustness of our algorithm to real-world operating conditions.

\section{Conclusion}
\label{sec:conclusion}

We present an algorithm for low-observability distribution system voltage estimation based on constrained low-rank matrix completion.
This method can accurately estimate voltage phasors under low-observability conditions where standard state estimation methods cannot operate, and can flexibly accommodate any distribution network measurements available in the field.
Our empirical evaluations of this method demonstrate that it produces accurate and robust voltage phasor estimates on the IEEE 33- and 123-bus test systems under a variety of data availability conditions.
As such, we believe that our algorithm is a useful mechanism for voltage estimation on modern distribution systems.

\bibliography{references}

\begin{thebibliography}{10}
\providecommand{\url}[1]{#1}
\csname url@samestyle\endcsname
\providecommand{\newblock}{\relax}
\providecommand{\bibinfo}[2]{#2}
\providecommand{\BIBentrySTDinterwordspacing}{\spaceskip=0pt\relax}
\providecommand{\BIBentryALTinterwordstretchfactor}{4}
\providecommand{\BIBentryALTinterwordspacing}{\spaceskip=\fontdimen2\font plus
\BIBentryALTinterwordstretchfactor\fontdimen3\font minus
  \fontdimen4\font\relax}
\providecommand{\BIBforeignlanguage}[2]{{%
\expandafter\ifx\csname l@#1\endcsname\relax
\typeout{** WARNING: IEEEtran.bst: No hyphenation pattern has been}%
\typeout{** loaded for the language `#1'. Using the pattern for}%
\typeout{** the default language instead.}%
\else
\language=\csname l@#1\endcsname
\fi
#2}}
\providecommand{\BIBdecl}{\relax}
\BIBdecl

\bibitem{liacco1982role}
T.~D. Liacco, ``The role of state estimation in power system operation,''
  \emph{IFAC Proceedings Volumes}, vol.~15, no.~4, pp. 1531--1533, 1982.

\bibitem{Abur2004}
\BIBentryALTinterwordspacing
A.~Abur and A.~G. Exposito, \emph{{Power System State Estimation: Theory and
  Implementation}}.\hskip 1em plus 0.5em minus 0.4em\relax Abingdon: Dekker,
  2004. [Online]. Available: \url{http://cds.cern.ch/record/994935}
\BIBentrySTDinterwordspacing

\bibitem{dehghanpour2018survey}
K.~Dehghanpour, Z.~Wang, J.~Wang, Y.~Yuan, and F.~Bu, ``A survey on state
  estimation techniques and challenges in smart distribution systems,''
  \emph{IEEE Transactions on Smart Grid}, 2018.

\bibitem{deng2002branch}
Y.~Deng, Y.~He, and B.~Zhang, ``A branch-estimation-based state estimation
  method for radial distribution systems,'' \emph{IEEE Transactions on power
  delivery}, vol.~17, no.~4, pp. 1057--1062, 2002.

\bibitem{pereira2004integrated}
J.~Pereira, J.~Saraiva, and V.~Miranda, ``An integrated load allocation/state
  estimation approach for distribution networks,'' in \emph{Probabilistic
  Methods Applied to Power Systems, 2004 International Conference on}.\hskip
  1em plus 0.5em minus 0.4em\relax IEEE, 2004, pp. 180--185.

\bibitem{driesen2008design}
J.~Driesen and F.~Katiraei, ``Design for distributed energy resources,''
  \emph{IEEE Power and Energy Magazine}, vol.~6, no.~3, 2008.

\bibitem{primadianto2017review}
A.~Primadianto and C.-N. Lu, ``A review on distribution system state
  estimation,'' \emph{IEEE Transactions on Power Systems}, vol.~32, no.~5, pp.
  3875--3883, 2017.

\bibitem{singh2009}
R.~Singh, B.~C. Pal, and R.~B. Vinter, ``Measurement placement in distribution
  system state estimation,'' \emph{IEEE Transactions on Power Systems},
  vol.~24, no.~2, pp. 668--675, 2009.

\bibitem{bhela2018}
S.~Bhela, V.~Kekatos, and S.~Veeramachaneni, ``Enhancing observability in
  distribution grids using smart meter data,'' \emph{IEEE Transactions on Smart
  Grid}, vol.~9, no.~6, pp. 5953--5961, 2018.

\bibitem{YCNN}
H.~Jiang and Y.~Zhang, ``Short-term distribution system state forecast based on
  optimal synchrophasor sensor placement and extreme learning machine,'' in
  \emph{2016 IEEE Power and Energy Society General Meeting (PESGM)}, July 2016,
  pp. 1--5.

\bibitem{manitsas2012}
E.~Manitsas, R.~Singh, B.~C. Pal, and G.~Strbac, ``Distribution system state
  estimation using an artificial neural network approach for pseudo measurement
  modeling,'' \emph{IEEE Transactions on Power Systems}, vol.~27, no.~4, pp.
  1888--1896, Nov 2012.

\bibitem{wu2013}
J.~Wu, Y.~He, and N.~Jenkins, ``A robust state estimator for medium voltage
  distribution networks,'' \emph{IEEE Transactions on Power Systems}, vol.~28,
  no.~2, pp. 1008--1016, May 2013.

\bibitem{clements}
K.~A. Clements, ``The impact of pseudo-measurements on state estimator
  accuracy,'' in \emph{2011 IEEE Power and Energy Society General Meeting},
  July 2011, pp. 1--4.

\bibitem{NNVoltages}
M.~Pertl, K.~Heussen, O.~Gehrke, and M.~Rezkalla, ``Voltage estimation in
  active distribution grids using neural networks,'' in \emph{2016 IEEE Power
  and Energy Society General Meeting (PESGM)}, July 2016, pp. 1--5.

\bibitem{ExactMatrix}
\BIBentryALTinterwordspacing
E.~J. Cand{\`e}s and B.~Recht, ``Exact matrix completion via convex
  optimization,'' \emph{Foundations of Computational Mathematics}, vol.~9,
  no.~6, pp. 717--772, Apr 2009. [Online]. Available:
  \url{https://doi.org/10.1007/s10208-009-9045-5}
\BIBentrySTDinterwordspacing

\bibitem{gao2016}
P.~{Gao}, M.~{Wang}, S.~G. {Ghiocel}, J.~H. {Chow}, B.~{Fardanesh}, and
  G.~{Stefopoulos}, ``Missing data recovery by exploiting low-dimensionality in
  power system synchrophasor measurements,'' \emph{IEEE Transactions on Power
  Systems}, vol.~31, no.~2, pp. 1006--1013, March 2016.

\bibitem{genes2018}
C.~{Genes}, I.~{Esnaola}, S.~M. {Perlaza}, L.~F. {Ochoa}, and D.~{Coca},
  ``Robust recovery of missing data in electricity distribution systems,''
  \emph{IEEE Transactions on Smart Grid}, pp. 1--1, 2018.

\bibitem{klauber}
C.~Klauber and H.~Zhu, ``Distribution system state estimation using
  semidefinite programming,'' in \emph{2015 North American Power Symposium
  (NAPS)}, Oct 2015, pp. 1--6.

\bibitem{TruncNorm}
Y.~Hu, D.~Zhang, J.~Ye, X.~Li, and X.~He, ``Fast and accurate matrix completion
  via truncated nuclear norm regularization,'' \emph{IEEE Transactions on
  Pattern Analysis and Machine Intelligence}, vol.~35, no.~9, pp. 2117--2130,
  Sept 2013.

\bibitem{LargeData}
\BIBentryALTinterwordspacing
R.~H. Keshavan, S.~Oh, and A.~Montanari, ``Matrix completion from a few
  entries,'' \emph{CoRR}, vol. abs/0901.3150, 2009. [Online]. Available:
  \url{http://arxiv.org/abs/0901.3150}
\BIBentrySTDinterwordspacing

\bibitem{AnyBasis}
\BIBentryALTinterwordspacing
D.~Gross, ``Recovering low-rank matrices from few coefficients in any basis,''
  \emph{CoRR}, vol. abs/0910.1879, 2009. [Online]. Available:
  \url{http://arxiv.org/abs/0910.1879}
\BIBentrySTDinterwordspacing

\bibitem{candes2010matrix}
E.~J. Candes and Y.~Plan, ``Matrix completion with noise,'' \emph{Proceedings
  of the IEEE}, vol.~98, no.~6, pp. 925--936, 2010.

\bibitem{Bernstein_Dall’Anese_2017}
\BIBentryALTinterwordspacing
A.~Bernstein and E.~Dall'Anese, ``Linear power-flow models in multiphase
  distribution networks: Preprint,'' May 2017. [Online]. Available:
  \url{http://www.osti.gov/scitech/servlets/purl/1361015}
\BIBentrySTDinterwordspacing

\bibitem{Liu2017}
Y.~Liu, N.~Zhang, Y.~Wang, J.~Yang, and C.~Kang, ``Data-driven power flow
  linearization: A regression approach,'' 2017,
  https://arxiv.org/pdf/1710.10621.

\bibitem{bernstein2017}
A.~Bernstein, C.~Wang, E.~Dall'Anese, J.-Y. {Le Boudec}, and C.~Zhao,
  ``Load-flow in multiphase distribution networks: Existence, uniqueness, and
  linear models,'' 2017, http://arxiv.org/abs/1702.03310.

\end{thebibliography}
\bibliographystyle{IEEEtran}

\ifCLASSOPTIONcaptionsoff
  \newpage
\fi




\vspace{-0.5cm}

\vfill
\newpage







\excludeforsubmit{
\input{appendices.tex}}

\end{document}